\documentclass[a4paper,12pt,english]{smfart}

\usepackage{amssymb}
\usepackage{url}
\usepackage[T1]{fontenc}
\RequirePackage{calrsfs}
\DeclareSymbolFont{rsfscript}{OMS}{rsfs}{m}{b}
\DeclareSymbolFontAlphabet{\mathrsfs}{rsfscript}
\usepackage[utopia,expert]{mathdesign}
\usepackage{lscape}
\usepackage{fancybox}

\usepackage[leftbars]{changebar}

\usepackage{palatino}
\usepackage{rotating}
\usepackage{graphicx}

\usepackage{enumerate}
\usepackage{enumitem}

\usepackage{color}
\definecolor{shadecolor}{gray}{0.90}

\def\bfit{\bfseries\itshape}

\def\proofname{D\'emonstration}

\makeindex

\def\indexname{Index des notations}

\newtheorem{theo}{Theorem}[section]
\def\thetheo{\thesection.\arabic{theo}}
\newtheorem{prop}[theo]{Proposition}

\newtheorem{lem}[theo]{Lemma}

\newtheorem{coro}[theo]{Corollary}

\newtheorem{defi}[theo]{Definition}

\renewcommand\theequation{\thetheo}
\def\equat{\refstepcounter{theo}\begin{equation}}
\def\endequat{\end{equation}}

\renewcommand\thetable{\thetheo}

\renewcommand\thesection{\arabic{section}}
\renewcommand\thesubsection{\thesection.\Alph{subsection}}

\def\contentsname{Table des mati\`eres}
\def\listfigurename{Liste des figures}
\def\listtablename{Liste des tables}
\def\appendixname{Appendice}

\def\refname{R\'ef\'erences}

\def\IG{{\mathfrak I}}

    \def\NM{{\mathbb{N}}}

\def\SG{{\mathfrak S}}

    \def\ZM{{\mathbb{Z}}}

  \def\ab{{\mathbf a}}  \def\AC{{\mathcal{A}}}
    
  \def\cb{{\mathbf c}}  
  \def\db{{\mathbf d}}  \def\DC{{\mathcal{D}}}
    \def\EC{{\mathcal{E}}}
    
  \def\gb{{\mathbf g}}  
    \def\HC{{\mathcal{H}}}

    \def\LC{{\mathcal{L}}}
    \def\MC{{\mathcal{M}}}
    
    \def\OC{{\mathcal{O}}}
    \def\PC{{\mathcal{P}}}
    
    \def\RC{{\mathcal{R}}}

    \def\XC{{\mathcal{X}}}

\def\Crm{{\mathrm{C}}}

  \def\frm{{\mathrm{f}}}  
    
\def\Hrm{{\mathrm{H}}}

\def\Lrm{{\mathrm{L}}}

  \def\wti{{\tilde{w}}}

          \def\aov{{\overline{a}}}

\def\Mov{{\overline{M}}}          \def\mov{{\overline{m}}}

\def\b{\beta}
\def\g{\gamma}
\def\G{\Gamma}
\def\d{\delta}

\def\ph{\varphi}
\def\l{\lambda}
\def\L{\Lambda}
\def\m{\mu}
\def\o{\omega}
\def\O{\Omega}
\def\r{\rho}
\def\s{\sigma}

\def\th{\theta}

\def\t{\tau}

\def\z{\zeta}

\def\alpb{{\boldsymbol{\alpha}}}

\def\mub{{\boldsymbol{\mu}}}

\DeclareMathOperator{\Id}{{\mathrm{Id}}}

\DeclareMathOperator{\val}{{\mathrm{val}}}

\def\to{\rightarrow}
\def\longto{\longrightarrow}

\def\fonction#1#2#3#4#5{\begin{array}{rccc}
{#1} : & {#2} & \longto & {#3}  \\
& {#4} & \longmapsto & {#5} 
\end{array}}

\def\fonctio#1#2#3#4{\begin{array}{ccc}
{#1} & \longto & {#2} \\
{#3} & \longmapsto & {#4} 
\end{array}}

\def\DS{\displaystyle}

\def\finl{~$\blacksquare$}

\def\lexp#1#2{\kern\scriptspace\vphantom{#2}^{#1}\kern-\scriptspace#2}
\def\le{\hspace{0.1em}\mathop{\leqslant}\nolimits\hspace{0.1em}}
\def\ge{\hspace{0.1em}\mathop{\geqslant}\nolimits\hspace{0.1em}}

\mathchardef\inferieur="321E
\mathchardef\superieur="321F

\def\eqna{\begin{eqnarray*}}
\def\endeqna{\end{eqnarray*}}

\def\itemth#1{\item[${\mathrm{(#1)}}$]}

\catcode`\@=11
\long\def\@car#1#2\@nil{#1}
\long\def\@first#1#2{#1}
\long\def\@second#1#2{#2}
\long\def\ifempty#1{\expandafter\ifx\@car#1@\@nil @\@empty
  \expandafter\@first\else\expandafter\@second\fi}
\catcode`\@=12

\def\boitegrise#1#2{\begin{centerline}{\fcolorbox{black}{shadecolor}{~
    \begin{minipage}[t]{#2}{\vphantom{~}#1\vphantom{$A_{\DS{A_A}}$}}
            \end{minipage}~}}\end{centerline}\medskip}

\def\surto{\twoheadrightarrow}

\theoremstyle{remark}
\newtheorem{rema}[theo]{Remark}

\newtheorem{exemple}[theo]{Example}

\theoremstyle{plain}

\def\BIL{LR}
\def\GAUCHE{L}
\def\CAR{CAR}
\def\FAM{FAM}

\def\xyinj{\ar@{^{(}-&gt;}}
\def\xysur{\ar@{-&gt;&gt;}}

\def\gauche{{\mathrm{left}}}
\def\droite{{\mathrm{right}}}

\bigskip

\makeatletter
\def\hlinewd#1{%
\noalign{\ifnum0=`}\fi\hrule \@height #1 %
\futurelet\reserved@a\@xhline}
\makeatother

\newlength\epaisLigne

\def\dotcup{\hskip1mm\dot{\cup}\hskip1mm}

\newcommand{\longiso}{\stackrel{\sim}{\longrightarrow}}

\makeatletter
\def\hlinewd#1{%
\noalign{\ifnum0=`}\fi\hrule \@height #1 %
\futurelet\reserved@a\@xhline}
\makeatother

\usepackage{lscape}

\def\irr{{\mathrm{ir,f}}}
\def\cactus{{\mathrm{Cact}}}

\def\lel{\hspace{0.03em}\mathop{\leqslant_L}\nolimits\hspace{0.03em}}
\def\leli{\hspace{0.03em}\mathop{\leqslant_L^I}\nolimits\hspace{0.03em}}
\def\leri{\hspace{0.03em}\mathop{\leqslant_R^I}\nolimits\hspace{0.03em}}
\def\lelri{\hspace{0.03em}\mathop{\leqslant_{LR}^I}\nolimits\hspace{0.03em}}
\def\ler{\hspace{0.03em}\mathop{\leqslant_R}\nolimits\hspace{0.03em}}
\def\lelr{\hspace{0.03em}\mathop{\leqslant_{LR}}\nolimits\hspace{0.03em}}

\begin{document}

\baselineskip=16pt

\title{Cells and cacti}

\author{{\sc C\'edric Bonnaf\'e}}
\address{
Institut Montpelli\'erain Alexander Grothendieck (CNRS: UMR 5149), 
Universit\'e Montpellier 2,
Case Courrier 051,
Place Eug\`ene Bataillon,
34095 MONTPELLIER Cedex,
FRANCE} 

\makeatletter
\email{cedric.bonnafe@univ-montp2.fr}
\makeatother

%
%
\date{\today}
\thanks{The author is partly supported by the ANR (Project No ANR-12-JS01-0003-01 ACORT)}
\maketitle
\pagestyle{myheadings}

\markboth{\sc C. Bonnaf\'e}{\sc Cells and cacti}

\def\opp{{\mathrm{op}}}

\begin{abstract}
Let $(W,S)$ be a Coxeter system, let $\ph$ be a weight function on $S$ and let $\cactus_W$ 
denote the associated {\it cactus group}. Following an idea of I. Losev, we construct 
an action of $\cactus_W \times \cactus_W$ on $W$ which has nice properties with respect 
to the partition of $W$ into left, right or two-sided cells (under some hypothesis, 
which hold for instance if $\ph$ is constant).
It must be noticed that the action depends heavily on $\ph$. 
\end{abstract}

Let $(W,S)$ be a Coxeter system with $S$ finite and let $\ph$ be a positive {\it weight function} 
on $S$ as defined by Lusztig~\cite{lusztig hecke}. We denote by $\cactus_W$ the {\it Cactus group} 
associated with $W$, as defined for instance in~\cite{losev} (see also Section~\ref{sec:cactus}). 
In~\cite{losev}, I. Losev 
has constructed, whenever $W$ is a finite Weyl group and $\ph$ is constant, an action of 
$\cactus_W \times \cactus_W$ on $W$ which satisfies some good properties with respect to the partition 
of $W$ into cells. His construction is realized as the combinatorial shadow of wall-crossing functors 
on the category $\OC$. 

In~\cite[\S{5.1}]{losev}, I. Losev suggested that 
this action could be obtained without any reference to some category $\OC$, and thus extended to other 
types of Coxeter groups and general weight functions $\ph$, using some recent results of 
Lusztig~\cite{lusztig action}. This is the aim of this paper to show that Losev's idea works, 
by using slight extensions of results from~\cite{bonnafe geck} and assuming that some of Lusztig's 
Conjectures in~\cite[\S{14.2}]{lusztig hecke} hold, as in~\cite{lusztig action}.
Note that, if $\ph$ is constant, then these Conjectures hold, so this provides 
at least an action in the equal parameter case: if moreover $W$ is a Weyl group, 
this action coincides with the one constructed by Losev~\cite{losev}. 

Let us now state our main result. If $I \subset S$, we denote by $W_I$ the subgroup 
generated by $I$ and by $\ph_I$ the restriction of $\ph$ to $I$. 
If $C$ is a left (respectively right) cell, then $\HC^L[C]$ 
(respectively $\HC^R[C]$) denotes the associated left (respectively right) $\HC$-module 
and $c_w^L$ (respectively $c_w^R$) denotes the image of the Kazhdan-Lusztig 
basis element $C_w$ in this module (see~\S\ref{sub:cells}). Finally, we set $\mub_2=\{1,-1\}$. 

\bigskip

\noindent{\bfit Theorem.--- } 
{\it Assume that {\bfit Lusztig's Conjectures P1, P4, P8 and P9 in~\cite[\S{14}]{lusztig hecke} hold 
for all triples $(W_I,I,\ph_I)$ such that $W_I$ is finite}. 
Then there exists an action of $\cactus_W \times \cactus_W$ 
on the set $W$ such that, if we denote by $\t_\ph^L$ (respectively $\t_\ph^R$) the permutation of $W$ obtained 
through the action of $(\t,1) \in \cactus_W \times \cactus_W$ (respectively $(1,\t) \in \cactus_W \times \cactus_W$), 
then:
\begin{itemize}
\itemth{a} If $C$ is a left cell, then $\t_\ph^L(C)$ is also a left cell. Moreover, 
there exists a sign map $\eta_L^{\t,\ph} : W \to \mub_2$ such that the $A$-linear map 
$\HC^L[C] \longiso \HC^L[\t_\ph^L(C)]$, $c_w^L \mapsto \eta_{L,w}^{\t,\ph} c_{\t_\ph^L(w)}^L$ is 
an isomorphism of left $\HC$-modules.

\itemth{a'} If $C$ is a right cell, then $\t_\ph^R(C)$ is a also right cell. Moreover, 
there exists a sign map $\eta_R^{\t,\ph} : W \to \mub_2$ such that the $A$-linear map 
$\HC^R[C] \longiso \HC^R[\t_\ph^R(C)]$, $c_w^R \mapsto \eta_{R,w}^{\t,\ph} c_{\t_\ph^R(w)}^R$ is 
an isomorphism of right $\HC$-modules.

\itemth{b} If $w \in W$, then $\t_\ph^L(w) \sim_R w$ and $\t_\ph^R(w) \sim_L w$.
\end{itemize}}

\bigskip

\noindent{\bfit Commentary.--- } Lusztig~\cite[\S{14.2}]{lusztig hecke} proposed several Conjectures 
relating the so-called {\it Lusztig's $\ab$-function} and the partition of $W$ into cells. Throughout 
this paper, the expression {\it Lusztig's Conjecture Pi} will refer 
to~\cite[\S{14.2},~Conjecture~Pi]{lusztig hecke} (for $1 \le i \le 15$). 
For instance, they all hold if $\ph$ is constant~\cite[\S{15}]{lusztig hecke}.\finl

\bigskip

\noindent{\bfit Acknowledgements.--- } I wish to thank warmly I. Losev for sending me 
his first version of~\cite{losev}, and for the e-mails we have exchanged afterwards.

\bigskip
\def\invo{{\overline{\hphantom{A}\vphantom{a}}}}

\section{Notation}

\medskip

\boitegrise{{\bf Set-up.} 
{\it We fix a Coxeter system $(W,S)$, whose length function is denoted by $\ell : W \to \NM$.
We also fix a totally ordered abelian group $\AC$ and we denote by $A$ the group algebra 
$\ZM[\AC]$. We use an exponential notation for $A$:
$$A=\oplus_{a \in \AC} \ZM v^a\quad\text{where}\quad v^av^{a'}=v^{a+a'}\quad\text{for all 
$a$, $a' \in \AC$.}$$
If $a_0 \in \AC$, we write $\AC_{\leqslant a_0}=\{a \in \AC~|~a \le a_0\}$ and 
$A_{\leqslant a_0}=\oplus_{a \in \AC_{\leqslant a_0}} \ZM v^a$;
we define similary $A_{<a_0}$, $A_{\geqslant a_0}$, $A_{> a_0}$. We denote by 
$\invo : A \to A$ the involutive automorphism such that $\overline{v^a}=v^{-a}$ 
for all $a \in \AC$. Since $\AC$ is totally ordered, $A$ inherits two maps $\deg : A \to \AC \cup \{-\infty\}$ 
and $\val : A \to \AC \cup \{+\infty\}$ respectively called {\bfit degree} and {\bfit valuation}, and 
which are defined as usual.\\
\hphantom{AA} We also fix a {\bfit weight function} $\ph : S \to \AC_{>0}$ 
(that is, $\ph(s)=\ph(t)$ for all $s$, $t \in S$ which are conjugate in $W$) and, if 
$I \subset S$, we denote by $\ph_I : I \to \AC_{>0}$ the restriction of $\ph$.}}{0.75\textwidth}

\bigskip

\subsection{Cells}\label{sub:cells}
Let $\HC=\HC(W,S,\ph)$ denote the Iwahori-Hecke algebra associated with the triple $(W,S,\ph)$. This 
$A$-algebra is free as an $A$-module, with a standard basis denoted by $(T_w)_{w \in W}$. 
The multiplication is completely determined by the following two rules:
$$
\begin{cases}
T_w T_{w'} = T_{ww'} &\text{if $\ell(ww')=\ell(w)+\ell(w')$,}\\
(T_s-v^{\ph(s)})(T_s + v^{-\ph(s)})=0 & \text{if $s \in S$.}\\
\end{cases}
$$The involution $\invo$ on $A$ can be extended to an $A$-semilinear involutive automorphism 
$\invo : \HC \to \HC$ by setting $\overline{T}_w=T_{w^{-1}}^{-1}$. 
Let 
$$\HC_{<0} = \mathop{\oplus}_{w \in W} A_{<0} T_w.$$
If $w \in W$, there exists~\cite{lusztig hecke} a unique $C_w \in \HC$ such that
$$
\begin{cases}
\overline{C}_w=C_w,\\
C_w \equiv T_w \mod \HC_{< 0}.
\end{cases}
$$
It is well-known~\cite{lusztig hecke} that $(C_w)_{w \in W}$ is an $A$-basis of $\HC$ (called 
the {\it Kazhdan-Lusztig basis}) and we will 
denote by $h_{x,y,z} \in A$ the structure constants, defined by
$$C_xC_y=\sum_{z \in W} h_{x,y,z} C_z.$$
We also write
$$C_y=\sum_{x \in W} p_{x,y}^* T_x,$$
with $p_{x,y}^* \in A$. Recall that $p_{y,y}^*=1$ and $p_{x,y}^* \in A_{<0}$ if $x \neq y$.

We will denote by $\lel$, $\ler$, $\lelr$, $<_L$, $<_R$, $<_{LR}$, $\sim_L$, $\sim_R$ and $\sim_{LR}$ 
the relations defined in~\cite{lusztig hecke} and associated with the triple $(W,S,\ph)$: 
the relation $\lel$ is the finest preorder on $W$ such that, for any $w \in W$, 
$\oplus_{x \lel w} AC_x$ 
is a left ideal of $\HC$, while $\sim_L$ is the associated equivalence relation associated 
(the other relations are defined similarly, by replacing left ideal by right or two-sided 
ideal). 
Also, we will call left, right and two-sided cells the equivalence classes for the relations 
$\sim_L$, $\sim_R$ and $\sim_{LR}$ respectively. If $C$ is a left cell, we set 
$$\HC^{\lel C}=\mathop{\oplus}_{w \lel C} A~C_w,\quad \HC^{<_L C}=\mathop{\oplus}_{w <_L C} A~C_w\quad\text{and}\quad
\HC^L[C]=\HC^{\lel C}/\HC^{<_L C}.$$
These are left $\HC$-modules. If $w \in C$, we denote by $c_w^L$ the image of $C_w$ in the quotient $\HC^L[C]$ 
and $\HC^{\lel C}$ and $\HC^{<_L C}$ might be also denoted by $\HC^{\lel w}$ and 
$\HC^{<_L w}$ respectively: 
it is clear that $(c_w^L)_{w \in C}$ is an $A$-basis of $\HC^L[C]$. 
If $C$ is a right (respectively two-sided) cell, we define similarly $\HC^{\ler C}$, $\HC^{<_R C}$ 
and $\HC^R[C]$ (respectively $\HC^{\lelr C}$, $\HC^{<_{LR} C}$ and $\HC^{LR}[C]$), 
as well as $c_w^R$ (respectively $c_w^{LR}$).

\bigskip
\def\proj{{\mathrm{pr}}}

\subsection{Parabolic subgroups} 
We denote by $\PC(S)$ the set of subsets of $S$. If $I \subset S$, we denote by $W_I$ the standard 
parabolic subgroup generated by $I$ and by $X_I$ the set of elements $x \in W$ which have minimal length 
in $x W_I$. We also define $\proj_L^I : W \to W_I$ and $\proj_R^I : W \to W_I$ by the following formulas:
$$\forall~x \in X_I,~\forall~w \in W_I,\quad\proj_L^I(xw)=w\quad\text{and}\quad\proj_R^I(wx^{-1})=w.$$
If $\d : W_I \to W_I$ is any map, we denote by $\d^L : W \to W$ and $\d^R : W \to W$ 
the maps defined by 
$$\d^L(xw)=x\d(w)\qquad\text{and}\qquad \d^R(wx^{-1})=\d(w)x^{-1}$$
for all $x \in X_I$ and $w \in W_I$ (see~\cite[\S{6}]{bonnafe geck}). 
We denote by $\d^\opp : W_I \to W_I$ the map defined by 
$$\d^\opp(w)=\d(w^{-1})^{-1}$$
for all $w \in W$. Note that $\d^R=((\d^\opp)^L)^\opp$. If $\s : W \to W$ is any automorphism such that 
$\s(S)=S$, then 
\equat\label{eq:auto-proj}
\s \circ \proj^I_L = \proj^{\s(I)}_L \circ \s\quad\text{and}\quad 
\s \circ \proj^I_R = \proj^{\s(I)}_R \circ \s.
\endequat
If $\EC$ is a set and $\mu : W_I \to \EC$ is any map, we define $\mu_L : W \to \EC$ 
(respectively $\mu_R : W \to \EC$) by 
$$\mu_L=\mu \circ \proj_L^I\qquad(\text{respectively~}\mu_R=\mu \circ \proj_R^I(w).$$
For instance, $\proj_L^I=(\Id_{W_I})_L$ and $\proj_R^I=(\Id_{W_I})_R$.

\medskip

The Hecke algebra $\HC(W_I,I,\ph_I)$ will be denoted by $\HC_I$ and will be viewed 
as a subalgebra of $\HC$ in the natural way. It follows from the multiplication rules in the 
Hecke algebra that the right $\HC_I$-module $\HC$ is free (hence flat) with basis $(T_x)_{x \in X_I}$. 
This remark has the following consequence (in the next lemma, if $E$ is a subset of $\HC$, 
then $\HC E$ denotes the left ideal generated by $E$):

\bigskip

\begin{lem}\label{lem:platitude}
If $\IG$ and $\IG'$ are left ideals of $\HC_I$ such that $\IG \subset \IG'$, then:
\begin{itemize}
\itemth{a} $\HC\IG=\oplus_{x \in X_I} T_x \IG$.

\itemth{b} The natural map $\HC \otimes_{\HC_I} \IG \to \HC\IG$ is an isomorphism of left $\HC$-modules.

\itemth{c} The natural map $\HC \otimes_{\HC_I} (\IG'/\IG) \to \HC \IG'/\HC \IG$ is an isomorphism 
of left $\HC$-modules.
\end{itemize}
\end{lem}

\bigskip

Let $\PC_\frm(S)$ (respectively $\PC_\irr(S)$) 
denote the set of subsets $I$ of $S$ such that $W_I$ is finite (respectively such that $W_I$ is finite 
and the Coxeter graph of $(W_I,I)$ is connected). If $I \in \PC_\frm(S)$, we denote 
by $w_I$ the longest element of $W_I$ and we set
$$\fonction{\o_I}{W_I}{W_I}{w}{w_Iww_I.}$$
It is an automorphism of $W_I$ which satisfies $\o_I(I) = I$. If $W$ is finite, then $w_S$ 
will be denoted by $w_0$, according to the tradition. Also, $\o_S$ will be denoted by $\o_0$. 

If $I \in \PC_\frm(S)$, we denote by $\ab_I : W_I \to \AC$ the {\it Lusztig's $\ab$-function} 
defined by
$$\ab_I(z) = \max_{x,y \in W_I} \deg(h_{x,y,z})$$
for all $z \in W_I$. We also set $\alpb_I(z)=\ab_I(w_I z)-\ab_I(z)$. 
If $W$ itself is finite, then $\ab_S$ and $\alpb_S$ will be simply denoted by $\ab$ and $\alpb$ 
respectively.

\bigskip

\subsection{Descent sets} 
If $w \in W$, we set
$$\LC(w)=\{s \in S~|~sw < w\}\qquad\text{and}\qquad
\RC(w)=\{s \in S~|~ws < w\}.$$
Then $\LC(w)$ (respectively $\RC(w)$) is called the {\it left descent set} 
(respectively {\it right descent set}) of $w$: it is easy to see that they both 
belong to $\PC_\frm(S)$. It is also well-known~\cite[Lemma~8.6]{lusztig hecke} that the map 
$\LC : W \to \PC_\frm(S)$ (respectively $\RC : W \to \PC_\frm(S)$) is constant 
on right (respectively left) cells.

\bigskip

\subsection{Cells and parabolic subgroups} 
We will now recall Geck's Theorem about the parabolic induction of cells~\cite{geck induction}. 
First, it is clear that $(C_w)_{w \in W_I}$ is the Kazhdan-Lusztig basis 
of $\HC_I$. We can then define a preorder $\leli$ and its associated equivalence class $\sim_L^I$ 
on $W_I$ in the same way as $\lel$ and $\sim_L$ are defined for $W$. We define similarly 
$\leri$, $\sim_R^I$, $\lelri$ and $\sim_{LR}^I$. 
If $w \in W$, then there exists a unique $a \in X_I$ and a unique 
$x \in W_I$ such that $w=ax$: we then set
$$G_w^I=T_a C_x.$$
It is easily seen that $(G_w^I)_{w \in W}$ is an $A$-basis of $\HC$ so that we can write, 
for $b \in X_I$ and $y \in W_I$,
$$C_{by}=\sum_{\substack{a \in X_I \\ x \in W_I}} p_{a,x,b,y}^I T_a C_x,$$
where $p_{a,x,b,y}^I \in A$. 


\bigskip

\begin{theo}[Geck]\label{theo:geck}
Let $E$ be a subset of $W_I$ such that, if $x \in E$ and if $y \in W_I$ is such that $y \leli x$, then 
$y \in E$. Let $\IG=\oplus_{w \in E} A~C_w$. Then
$$\HC \IG = \mathop{\oplus}_{w \in X_I \cdot E} A~G_w^I = \mathop{\oplus}_{w \in X_I \cdot E} A~C_w.$$
In particular, if $w$, $w'$ are elements of $W$ are such that $w \lel w'$ (respectively $w \sim_L w'$), 
then $\proj_L^I(w) \leli \proj_L^I(w')$ (respectively $\proj_L^I(w) \sim_L^I \proj_L^I(w')$).

Moreover, if $a$, $b \in X_I$ and $x$, $y \in W_I$, then:
\begin{itemize}
\itemth{a} $p_{b,y,b,y}^I=1$.

\itemth{b} If $ax \neq by$, then $p_{a,x,b,y}^I \in A_{<0}$.

\itemth{c} If $ax \neq by$ and $p_{a,x,b,y}^I \neq 0$, then $a < b$, $ax \le by$ and $x \leli y$.
\end{itemize}
\end{theo}

\bigskip

\begin{coro}[Geck]\label{coro:geck}
We have:
\begin{itemize}
\itemth{a} $\leli$ and $\sim_L^I$ are just the restriction of $\lel$ and $\sim_L$ to $W_I$ 
(and so we will use only the notation $\lel$ and $\sim_L$).

\itemth{b} If $C$ is a left cell in $W_I$, then $X_I \cdot C$ is a union of left cells of $W$.
\end{itemize}
\end{coro}

\section{Preliminaries}\label{section:invo}

\bigskip

\boitegrise{{\bf Hypothesis and notation.} {\it In this section, and only in this section, 
we fix an $A$-module $\MC$ and we assume that:
\begin{itemize}
\itemth{I1} $\MC$ admits an $A$-basis $(m_x)_{x \in X}$, where $X$ is a {\it poset}. 
We set 
$$\MC_{<0} = \oplus_{x \in X} A_{<0} m_x.$$
\itemth{I2} $\MC$ admits a semilinear involution $\overline{\hphantom{x}\vphantom{a}} : \MC \to \MC$. 
We set
$$\MC_{\mathrm{skew}} = \{m \in \MC~|~m+\overline{m}=0\}.$$
\itemth{I3} If $x \in X$, then $\DS{\overline{m}_x \equiv m_x 
\mod\Bigl(\mathop{\oplus}_{y < x} A m_y\Bigr)}$
\itemth{I4} If $x \in X$, then the set $\{y \in X~|~y \le x\}$ is finite.
\end{itemize}}\vskip-0.5cm}{0.75\textwidth}

\bigskip

\def\invotext{$\overline{\hphantom{x}\vphantom{a}}$}

\begin{prop}\label{prop:inf-anti}
The $\ZM$-linear map 
$$\fonctio{\MC_{<0}}{\MC_{\mathrm{skew}}}{m}{m-\overline{m}}$$
is an isomorphism.
\end{prop}

\bigskip

\begin{proof}
First, note that the corresponding result for the $A$-module $A$ itself holds. In other words, 
\equat\label{eq:inf-anti}
\text{\it The map $A_{<0} \to A_{\mathrm{skew}}$, $a \mapsto a -\overline{a}$ is an isomorphism.}
\endequat
Indeed, if $a \in A_{\mathrm{skew}}$, write $a=\sum_{\g \in \G} r_\g v^\g$, with $r_\g \in \ZM$. Now, 
if we set $a_-=\sum_{\g < 0} r_\g v^\g \in A_{<0}$, then $a=a_- - \overline{a}_-$. This shows the surjectivity, 
while the injectivity is trivial.

\medskip

Now, let $\L : \MC_{<0} \to \MC_{\mathrm{skew}}$, $m \mapsto m - \overline{m}$. 
For $\XC \subset X$, we set $\MC^\XC=\oplus_{x \in \XC} A \, m_x$ and $\MC_{<0}^\XC=\oplus_{x \in \XC} A_{<0}\, m_x$. 
Assume that, for all $x \in \XC$ and all $y \in X$ such that $y \le x$, then $y \in \XC$. By (I3), 
$\MC^\XC$ is stabilized by the involution ~\invotext. 
Since $X$ is the union of such finite $\XC$ (by (I4)), it shows that we may, and we will, 
assume that $X$ is finite. Let us write $X=\{x_0,x_1,\dots,x_n\}$ in such a way that, 
if $x_i \le x_j$, then $i \le j$ (this is always possible). For simplifying notation, 
we set $m_{x_i}=m_i$. Note that, by (I3), 
$$\overline{m}_i \in m_i + \Bigl(\mathop{\oplus}_{0 \le j < i} A\, m_j\Bigr).\leqno{(*)}$$
In particular, $\overline{m}_0 = m_0$. 

\medskip

Now, let $m \in \MC_{<0}$ be such that $\overline{m}=m$ and assume that $m \neq 0$. 
Write $m=\sum_{i=0}^r a_i m_i$, with $r \le n$, $a_i \in A_{<0}$ and $a_r \neq 0$. Then, by (I2), 
$$\overline{m} \equiv \overline{a}_r m_r \mod\Bigl(\mathop{\oplus}_{0 \le j <i} A m_j\Bigr).$$
Since $\overline{m}=m$, this forces $\overline{a}_r=a_r$, which is impossible 
(because $a_r \in A_{<0}$ and $a_r \neq 0$). So $\L$ is injective.

\medskip

Let us now show that $\L$ is surjective. So, 
let $m \in \MC_{\mathrm{skew}}$, and assume that $m \neq 0$ (for otherwise there is nothing to prove). 
Write $m=\sum_{i=0}^r a_i m_{i}$, with $r \le n$, $a_i \in A$ and $a_r \neq 0$. 
We shall prove by induction on $r$ that there exists $\mu \in \MC_{<0}$ such that $m=\mu - \overline{\mu}$. 
If $r=0$, then the result follows from~(\ref{eq:inf-anti}) and the fact that $\overline{m}_0=m_0$. 
So assume that $r > 0$. Then 
$$m + \overline{m} \equiv (a_r + \overline{a}_r) m_{r} \mod \MC^{\XC_{r-1}},$$
where $\XC_j=\{x_0,x_1,\dots,x_j\}$. 
Since $m + \overline{m}=0$, this forces $a_r \in A_{\mathrm{skew}}$. So, by~(\ref{eq:inf-anti}), 
there exists $a \in A_{<0}$ such that $a-\overline{a} = a_r$. 
Now, let $m'= m - a m_{r} + \overline{a} \overline{m}_{r}$. Then 
$m' + \overline{m}' = 0$ and $m' \in \oplus_{0 \le j < r} A \, m_j$. So, by the induction hypothesis, 
there exists $\mu' \in \MC_{<0}$ such that $m'=\mu' - \overline{\mu}'$. 
Now, set $\mu = a m_{r} + \mu'$. Then $\mu \in \MC_{<0}$ and $m=\mu-\overline{\mu}=\L(\mu)$, 
as desired.                                                   
\end{proof}

\bigskip

\begin{coro}\label{coro:anti-stable}
Let $m \in \MC$. Then there exists a unique $M \in \MC$ such that 
$$
\begin{cases}
\Mov=M,\\
M \equiv m \mod \MC_{< 0}.\\
\end{cases}
$$
\end{coro}

\bigskip

\begin{proof}
Setting $M=m+\mu$, the problem is equivalent to find $\mu \in \MC_{<0}$ such that $\overline{m+\mu} = m+\mu$. 
This is equivalent to find $\mu \in \MC_{< 0}$ such that 
$\mu-\overline{\mu} = \mov-m$: since $\mov - m \in \MC_{\mathrm{skew}}$, this problem 
admits a unique solution, thanks to Proposition~\ref{prop:inf-anti}. 
\end{proof}

\bigskip

The Corollary~\ref{coro:anti-stable} can be applied to the $A$-module 
$A$ itself. However, in this case, its proof becomes obvious: 
if $a_\circ = \sum_{\g \in \G} a_\g v^\g$, 
then $a=\sum_{\g \le 0} a_\g v^\g + \sum_{\g > 0} a_{-\g} v^\g$ is the unique 
element of $A$ such that $\aov=a$ and $a \equiv a_\circ \mod A_{<0}$.


\bigskip

\begin{coro}\label{coro:kl-base-triangulaire}
Let $\XC$ be a subset of $X$ such that, if $x \le y$ and $y \in \XC$, then $x \in \XC$. 
Let $M \in \MC$ be such that $\Mov=M$ and $M \in \MC^\XC + \MC_{<0}$. Then 
$M \in \MC^\XC$. 
\end{coro}

\bigskip

\begin{proof}
Let $M_0 \in \MC^\XC$ be such that $M \equiv M_0 \mod \MC_{<0}$. 
From the existence statement of Corollary~\ref{coro:anti-stable} applied to $\MC^\XC$, there 
exists $M' \in \MC^\XC$ such that $\Mov'=M'$ and 
$M' \equiv M_0 \mod \MC_{<0}^\XC$. The fact that $M=M' \in \MC^\XC$ 
now follows from the uniquenes statement of Corollary~\ref{coro:anti-stable}.
\end{proof}

\bigskip

\begin{coro}\label{coro:kl-general}
Let $x \in X$. Then there exists a unique element $M_x \in \MC$ 
such that 
$$\begin{cases}
\overline{M}_x = M_x,\\
M_x \equiv m_x \mod \MC_{<0}.
\end{cases}$$
Moreover, $M_x \equiv m_x  \mod \mathop{\oplus}_{y < x} A_{<0} m_y$ and 
$(M_x)_{x \in X}$ is an $A$-basis of $\MC$.
\end{coro}

\bigskip

\begin{proof}
The existence and uniqueness of $M_x$ follow from Corollary~\ref{coro:anti-stable}. 
The statement about the base change follows by applying this existence and uniqueness 
to $\MC^{X_x}$, where $X_x=\{y \in X~|~ y \le x\}$.

Finally, the fact that $(M_x)_{x \in X}$ is an $A$-basis of $\MC$ follows from the fact that the base 
change from $(m_x)_{x \in X}$ to $(M_x)_{x \in X}$ is unitriangular.
\end{proof}

\bigskip

\section{Cellular pairs} 

\medskip

We set $\mub_2=\{1,-1\}$. 
The following definition extends slightly~\cite[Definition~4.1]{bonnafe geck}:

\bigskip

\begin{defi}\label{defi:cellular}
Let $\d : W \to W$ and $\mu : W \to \mub_2$, $w \mapsto \mu_w$ be two maps. 
Then the pair $(\d,\mu)$ 
is called {\bfit left cellular} 
if the following conditions are satisfied for every left cell $C$ of $W$:
\begin{itemize}[leftmargin=1.3cm]
\itemth{LC1} $\d(C)$ is also is a left cell.

\itemth{LC2} The $A$-linear map $(\d,\mu)_C : \HC^L[C] \to \HC^L[\d(C)]$, 
$c_w^L \mapsto \mu_w c_{\d(w)}^L$ is an isomorphism of left $\HC$-modules.
\end{itemize}
It is called {\bfit strongly left cellular} if it is left cellular and if satisfies moreover 
the following condition:
\begin{itemize}[leftmargin=1.3cm]
\itemth{LC3} If $w \in W$, then $\d(w) \sim_R w$.
\end{itemize}
If $\mu$ is constant and $\d$ satisfies $(\Lrm\Crm 1)$ and $(\Lrm\Crm 2)$ (respectively 
$(\Lrm\Crm 1)$, $(\Lrm\Crm 2)$ and $(\Lrm\Crm 3)$), 
then we say that $\d$ is a {\bfit left cellular} map (respectively a {\bfit strongly 
left cellular} map). 

We define similarly the notions of {\bfit right cellular} and {\bfit strongly right cellular} 
pair or map, as well as the notion of {\bfit two-sided cellular} pair or map. 
\end{defi}

\bigskip

The case where $\mu$ is constant corresponds to~\cite[Definition~4.1]{bonnafe geck}.
We will see in the next section that there exist left cellular pairs $(\d,\mu)$ such that 
$\m$ is not constant.

\medskip

\subsection{Strongness} 
It is unclear if there exist left cellular pairs or maps which are not strongly left cellular. 
At least, we are able to show that this probably cannot happen in finite Coxeter groups:

\bigskip

\begin{prop}\label{prop:strongly}
Assume that $W$ is {\bfit finite} and that {\bfit Lusztig's Conjectures P4 and P9 hold 
for $(W,S,\ph)$}. Then any left (respectively right) cellular pair is 
strongly left (respectively right) cellular.
\end{prop}

\bigskip

\begin{proof}
Assume that $W$ is finite. Let $(\d,\mu)$ be a left cellular pair and let $C$ be a left cell of $W$. 
Let $K$ denote the fraction field of $A$. 
Since the algebra $K\HC=K \otimes_A \HC$ 
is semisimple, there exist two idempotents $e$ and $f$ of $K\HC$ such that 
$$K\HC^{\lel C} = K\HC e \oplus K\HC^{<_L C}
\qquad 
\text{and}\qquad 
K\HC^{\lel \d(C)} = K\HC f \oplus K\HC^{<_L \d(C)}.$$
If $w \in C$ (respectively $w \in \d(C)$), we write $C_w = c_w^e + d_w^e$ 
(respectively $C_w=c_w^f + d_w^f$) where $c_w^e \in K\HC e$ and $d_w^e \in K\HC^{<_L C}$ 
(respectively $c_w^f \in K\HC f$ and $d_w^f \in K\HC^{<_L \d(C)}$).  
Then, by hypothesis, the $K$-linear map $\d^* : K\HC e \longiso K\HC f$ such that 
$\d^*(c_w^e)=\mu_w c_{\d(w)}^f$ for all $w \in C$ is an isomorphism of $K\HC$-modules. 

Recall that any morphism of left $K\HC$-modules $K\HC e \to K\HC f$ is of the form $m \mapsto mh$ 
for some $h \in eK\HC f$. So 
there exists $h \in eK\HC f$ such that, 
for all $w \in C$, $c_w^e h = \mu_w c_{\d(w)}^f$. In other words,
$$C_w h - \mu_w C_{\d(w)} = d_w^e h - \mu_w d_{\d(w)}^f.$$
Now, let $\G$ denote the two-sided cell containing $C$. By the semisimplicity of $K\HC$ and the fact 
that $\HC^L[C] \simeq \HC^L[\d(C)]$, this forces $\d(C)$ to be contained in $\G$. 
By P4 and P9, we then have $d_w^e$, $d_{\d(w)}^f \in K\HC^{<_{LR}\G}$, and so
$$C_w h - \mu_w C_{\d(w)} \in K\HC^{<_{LR}\G}.$$
In particular, $\d(w) \ler w$. Similarly, $w \ler \d(w)$ and so $\d(w) \sim_R w$, as desired.
\end{proof}

\bigskip

Note also the following result:

\bigskip

\begin{prop}\label{prop:left-descent}
Let $(\d,\mu)$ be a left (respectively right) cellular pair and let $w \in W$. 
Then $\LC(\d(w))=\LC(w)$ (respectively $\RC(\d(w))=\RC(w)$). 
\end{prop}

\bigskip

\begin{proof}
Let $C$ denote the left cell of $w$ and let $s \in S$. Then $s \in \LC(w)$ if and only 
if $C_s c_w^L=(v^{\ph(s)} + v^{-\ph(s)}) c_w^L$. So the result follows from the fact 
that the map $(\d,\mu)_C$ is an isomorphism of left $\HC$-modules. 
\end{proof}

\bigskip

\subsection{Induction of cellular pairs} 
The next result extends slightly~\cite[Theorem~6.2]{bonnafe geck}. We present here 
a somewhat different proof, based on the results of Section~\ref{section:invo}. 

\bigskip

\begin{theo}\label{theo:induction-left-cellular}
Let $I$ be a subset of $S$ and let $(\d,\mu)$ be a left cellular pair for $(W_I,I,\ph_I)$. Then 
$(\d^L,\mu_L)$ is a left cellular pair for $(W,S,\ph)$. If moreover $(\d,\mu)$ is strongly 
left cellular, then $(\d^L,\mu_L)$ is strongly left cellular. 
\end{theo}

\bigskip

\begin{proof}
The proof is divided in several steps:

\medskip

\noindent{\it $\bullet$ First step: construction and properties of an isomorphism of left $\HC$-modules.} 
Let $C$ be a left cell of $W_I$. We denote by $\EC$ (respectively $\EC^\#$) 
the set of elements $w$ in $W_I$ such that $w \lel C$ (respectively $w <_L C$). 
By Lemma~\ref{lem:platitude} and  Theorem~\ref{theo:geck}, the families $(G_w^I)_{w \in X_I \cdot \EC}$ and 
$(C_w)_{w \in X_I \cdot \EC}$ are $A$-basis of $\HC \HC_I^{\lel C}$. Similarly, 
the families $(G_w^I)_{w \in X_I \cdot \EC^\#}$ and 
$(C_w)_{w \in X_I \cdot \EC^\#}$ are $A$-basis of $\HC \HC_I^{<_L C}$. 

If $w \in X_I \cdot C$, we denote by $\gb_w^I$ (respectively $\cb_w^I$) the image of 
$G_w^I$ (respectively $C_w$) in $\HC \HC_I^{\lel C} /\HC \HC_I^{<_L C}$. Again by 
Lemma~\ref{lem:platitude}, 
$$\HC \otimes_{\HC_I} \HC_I^L[C] \simeq \HC \HC_I^{\lel C} /\HC \HC_I^{<_L C}.$$
Therefore, $(\gb_w^I)_{w \in X_I \cdot C}$ and $(\cb_w^I)_{w \in X_I \cdot C}$ 
can be viewed as $A$-bases of $\HC \otimes_{\HC_I} \HC_I^L[C]$. 

Since the pair $(\d,\mu)$ is left cellular, the $A$-linear map 
$\HC_I^L[C] \to \HC_I^L[\d(C)]$, $c_w^L \mapsto \mu_w c_{\d(w)}^L$ is an isomorphism 
of left $\HC_I$-modules. Therefore, the $A$-linear map 
$$\fonction{\th}{\HC \otimes_{\HC_I} \HC_I^L[C]}{\HC \otimes_{\HC_I} \HC_I^L[\d(C)]}{\gb_w^I}{\mu_{L,w} \gb_{\d(w)}^I}$$
is an isomorphism of left $\HC$-modules. 

Now, the left $\HC$-modules $\HC \HC_I^{\lel C}$ and $\HC \HC_I^{<_L C}$ are stable under 
the involution~\invotext. So $\HC \otimes_{\HC_I} \HC_I^L[C]$ inherits an action of the 
involution~\invotext. Similarly, $\HC \otimes_{\HC_I} \HC_I^L[\d(C)]$ inherits an action of the 
involution~\invotext. Moreover, these two $A$-modules (endowed with~\invotext) satisfy the hypotheses 
(I1), (I2), (I3) and (I4) of Section~\ref{section:invo} (by Theorem~\ref{theo:geck}).

Also, it follows from the definition that the isomorphism $\th$ 
commutes with this involution. Therefore, $\overline{\th(\cb_w^I)}=\th(\cb_w^I)$ 
for all $w \in X_I \cdot C$. Moreover, it follows from Theorem~\ref{theo:geck} that 
$$\th(\cb_w^I) \equiv \mu_{L,w} \gb_{\d(w)}^I \mod \oplus_{x \in X_I \cdot \d(C)} A_{<0} \gb_x^I.$$
But the element $\cb_{\d(w)}^I$ is stable under the involution~\invotext~and, 
again by Theorem~\ref{theo:geck}, it satisfies
$$\cb_{\d(w)}^I \equiv \gb_{\d(w)}^I \mod \oplus_{x \in X_I \cdot \d(C)} A_{<0} \gb_x^I.$$
Therefore, by Proposition~\ref{prop:inf-anti},
\equat\label{eq:theta}
\th(\cb_w^I) = \mu_{L,w} \cb_{\d(w)}^I.
\endequat

\medskip

\noindent{\it $\bullet$ Second step: partition into left cells.} 
Now, assume that $w \sim_L w'$. 
According to Corollary~\ref{coro:geck}(b), there exists a unique 
cell $C$ in $W_I$ such that $w$, $w' \in X_I\cdot C$. 
By the definition of $\lel$ and $\sim_L$, there exist four sequences $x_1$,\dots, $x_m$, $y_1$,\dots, $y_n$, 
$w_1$,\dots, $w_m$, $w_1'$,\dots, $w_n'$ such that:
$$
\begin{cases}
w_1=w, w_m=w',\\
w_1'=w',w_n'=w,\\
\forall~i \in \{1,2,\dots,m-1\},~h_{x_i,w_i,w_{i+1}} \neq 0,\\
\forall~j \in \{1,2,\dots,n-1\},~h_{y_j,w_j',w_{j+1}'} \neq 0.\\
\end{cases}
$$
Therefore, we have $w'=w_m \lel \cdots \lel w_2 \lel w_1 = w =w_n' \lel \cdots \lel w_2' \lel w_1'=w$ 
and so 
$w=w_1 \sim_L w_2 \sim_L \cdots \sim_L w_m=w'=w_1' \sim_L w_2' \sim_L \cdots \sim_L w_n'=w$. 
Again by Corollary~\ref{coro:geck}(b), $w_i$, $w_j' \in X_I\cdot C$. So it follows 
from~(\ref{eq:theta}) that $h_{x,\d^L(w_i),\d^L(w_{i+1})}=\mu_{L,w_i}\mu_{L,w_{i+1}} h_{x,w_i,w_{i+1}}$ 
and $h_{x,\d^L(w_j'),\d^L(w_{j+1}')}=\mu_{L,w_j'}\mu_{L,w_{j+1}'} h_{y_j,w_j',w_{j+1}'}$ for all $x \in W$. 
Therefore,
$$
\begin{cases}
\forall~i \in \{1,2,\dots,m-1\},~h_{x_i,\d^L(w_i),\d^L(w_{i+1})} \neq 0,\\
\forall~j \in \{1,2,\dots,n-1\},~h_{y_j,\d^L(w_j'),\d^L(w_{j+1}')} \neq 0.\\
\end{cases}
$$
It then follows that 
\eqna
\d^L(w')=\d^L(w_m) \lel \cdots \lel \d^L(w_2) 
&\lel\!\!\! & \d^L(w_1)=\d^L(w)=\d^L(w_n')  \\
&\lel\!\!\!& \cdots \lel \d^L(w_2') \lel \d^L(w_1')=\d^L(w'),
\endeqna
and so $\d^L(w) \sim_L \d^L(w')$, as expected. So we have proved that 
$$\text{\it if $w \sim_L w'$, then $\d^L(w) \sim_L \d^L(w')$.}\leqno{(*)}$$
Now, let $\d_1 : W_I \to W_I$ be the map defined by $\d_1(x)=x$ if $x \not\in \d(C)$ 
and $\d_1(\d(x))=x$ if $x \in C$. Let $\mu_1 : W \to \mub_2$ be defined by $\mu_{1,x}=1$ if 
$x \not\in \d(C)$ and $\mu_{1,\d(x)}=\mu_x$ if $x \in C$. Since left cellular maps 
can be defined ``locally'' (i.e. left cells by left cells), it is easily checked that 
$(\d_1,\mu_1)$ is left cellular. So, applying~$(*)$ to the pair $(\d_1,\mu_1)$ 
with $w$ and $w'$ replaced by $\d^L(w)$ and $\d^L(w')$, we obtain 
\equat\label{eq:left-cellular-bon}
\text{\it $w \sim_L w'$ if and only if $\d^L(w) \sim_L \d^L(w')$.}
\endequat

\medskip

\noindent{$\bullet$ \it Third step: left cellularity.} 
Now, let $C'$ be a left cell in $W$. It follows from~(\ref{eq:left-cellular-bon}) that 
$\d^L(C')$ is also a left cell and it follows from~(\ref{eq:theta}) that the $A$-linear 
map $\HC^L[C'] \to \HC^L[\d(C')]$, $c_w^L \mapsto \mu_{L,w} c_{\d^L(w)}^L$ is an isomorphism 
of left $\HC$-modules. In other words, $(\d^L,\mu_L)$ is left cellular.

\medskip

\noindent{\it $\bullet$ Fourth step: strongness.} 
Assume moreover that $(\d,\mu)$ is strongly left cellular. Let $w \in W$. Let us write 
$w=ax$ with $a \in X_I$ and $x \in W_I$. Then $\d(x) \sim_R x$ by~(LC3) 
and so $\d^L(w)=a\d(x) \sim_R a x=w$ by~\cite[Proposition~9.11]{lusztig hecke}.
\end{proof}

\bigskip

The next result extends slightly~\cite[Lemma~3.8]{geck relative}.

\begin{coro}\label{coro:egalite-geck}
Let $(\d,\mu)$ be a left cellular pair for $(W_I,I,\ph_I)$ and let $a$, $b \in X_I$ and 
$x$, $y \in W_I$ be such that $x \sim_L y$. Then
$$p_{a,x,b,y}^I=\mu_x\mu_y p_{a,\d(x),b,\d(y)}^I.$$
\end{coro}

\bigskip

\begin{proof}
This follows from~(\ref{eq:theta}).
\end{proof}

\bigskip

\section{Action of the longest element}

\medskip

\boitegrise{{\bf Hypothesis.} {\it We fix in this section a subset $I \in \PC_\frm(S)$ 
such that {\bfit Lusztig's Conjectures~P1, P4,~P8 and~P9 hold} 
for the triple $(W_I,I,\ph_I)$.}}{0.75\textwidth}

\bigskip

\begin{exemple}\label{rem:equal-bon}
Recall from~\cite[\S{15}]{lusztig hecke} that, if the weight function $\ph_I$ is constant, then 
Lusztig's Conjectures~P1,~P2,~P3,\dots,~P15 hold for $(W_I,I,\ph_I)$.\finl
\end{exemple}

\bigskip

\subsection{} 
The following result (which is crucial for our purpose) has been proved by Mathas~\cite{mathas} 
in the equal parameter case and extended by Lusztig~\cite[Theorem~2.3]{lusztig action} 
in the unequal parameter case:

\bigskip

\begin{theo}[Mathas, Lusztig]\label{theo:lusztig}
Let $I \in \PC_\frm(S)$ be such that {\bfit Lusztig's Conjectures P1, P4,~P8 and~P9 hold} 
for the triple $(W_I,I,\ph_I)$. Then there exists a (unique) sign map $\eta^I : W_I \to \mub_2$, $w \mapsto \eta_w^I$ 
and two (unique) involutions $\r_I$ and $\l_I$ 
of the set $W_I$ such that, for all $w \in W_I$, 
$$v^{\alpb_I(w)} T_{w_I} C_w \equiv \eta_w^I C_{\r_I(w)} \mod \HC_I^{<_{LR}^I w}$$
$$v^{\alpb_I(w)} C_w T_{w_I} \equiv \eta_w^I C_{\l_I(w)} \mod \HC_I^{<_{LR}^I w}.\leqno{\text{\it and}}$$
Note that $\l_I=\r_I^\opp$, that $\r_I=\l_I \circ \o_I$ and that
$$\r_I(w) \sim_L w\qquad\text{and}\qquad\l_I(w) \sim_R w.$$
\end{theo}

\bigskip

If $W$ itself is finite and if Lusztig's Conjectures~P1,~P4,~P8 and~P9 hold for $(W,S,\ph)$, 
then $\l_S$, $\r_S$ and $\eta^S$ will simply be denoted by $\l$, $\r$ and $\eta$ respectively. 

\bigskip

\begin{rema}\label{rem:simplification}
We will explain here why we only need to assume that Lusztig's Conjectures~P1,~P4,~P8 and~P9 hold 
for the above Theorem to hold (in~\cite[Theorem~2.5]{lusztig action}, Lusztig assumed that~P1,~P2,\dots,~P14 
and~P15 hold). This will be a consequence of a simplification of the proof of~\cite[Lemma~1.13]{lusztig action}, 
based on the ideas of~\cite{bonnafe two}. In particular, we avoid the use of the difficult 
Lusztig's Conjecture~P15 and the construction/properties of the asymptotic algebra. 

So assume that Lusztig's Conjectures~P1,~P4,~P8 and~P9 hold. We may, and we will, assume that 
$I=S$ (for simplifying notation). Let us write
$$T_{w_0} C_y = \sum_{x \lel y} \l_{x,y}~C_x,$$
with $\l_{x,y} \in A$. Note that
$$T_{w_0}^{-1} C_y = \sum_{x \lel y} \overline{\l}_{x,y}~C_x.$$
By~\cite[Proposition~1.4(a)]{bonnafe two}, 
$$\text{$\deg(\l_{x,y}) \le -\alpb(x)$ with equality only if $x \sim_L y$.}$$
By~\cite[Proposition~1.4(b)]{bonnafe two}, 
$$\text{$\deg(\overline{\l}_{x,y}) \le \alpb(y)$ with equality only if $x \sim_L y$.}$$
Assume now that $x \sim_L y$. Then $\alpb(x)=\alpb(y)$ by~P4 
and~\cite[Corollary~11.7]{lusztig hecke}, so 
$$\deg(\l_{x,y}) \le -\alpb(y) \le \val(\l_{x,y}).$$
So 
$$\text{if $x \sim_L y$, then $v^{\alpb(y)} \l_{x,y} \in \ZM$,}$$
Thanks to~P9, this is exactly the statement in~\cite[Lemma~1.13(a)]{lusztig action}. Note also 
that~\cite[Lemma~1.13(b)]{lusztig action} is already proved in~\cite[Proposition~1.4(c)]{bonnafe two}. 

One can then check that, once~\cite[Lemma~1.13]{lusztig action} is proved, 
the argument developed in~\cite[Proof~of~Theorem~2.3]{lusztig action} 
to obtain Theorem~\ref{theo:lusztig} does not make use any more 
of Lusztig's Conjectures.\finl
\end{rema}

\bigskip

\begin{rema}\label{rem:signe-constant}
In the equal parameter case, Mathas proved moreover that the sign map $w \mapsto \eta_w^I$ is constant on 
two-sided cells. However, this property does not hold in general, as it can be seen from direct 
computations whenever $W$ is of type $B_3$ (and $\ph$ is given by $\ph(t)=2$ and $\ph(s_1)=\ph(s_2)=1$, 
where $S=\{t,s_1,s_2\}$ and $s_1s_2$ has order $3$).\finl
\end{rema}

\bigskip

\begin{exemple}\label{ex:1}
Assume here that $W$ is {\it finite}. 
Since $\{1\}$ and $\{w_0\}$ are two-sided cells, we have $\l(1)=\r(1)=1$ and $\l(w_0)=\r(w_0)=w_0$. 
Moreover, $\eta_1=(-1)^{\ell(w_0)}$ and $\eta_{w_0}=1$.\finl
\end{exemple}

\bigskip

\subsection{Cellularity}
One of the key results towards a construction of an action of the cactus group 
is the following: 

%

\bigskip

\begin{theo}\label{theo:left-cellular}
Let $I \in \PC_\frm(S)$ be such that {\bfit Lusztig's Conjectures~P1, P4,~P8 and~P9 hold} 
for the triple $(W_I,I,\ph_I)$. Then the pair $(\l_I,\eta^I)$ (respectively $(\r_I,\eta^I)$) 
is strongly left (respectively right) cellular.
\end{theo}

\bigskip

\begin{proof}
For simplifying notation, we may, and we will, assume that $W$ is finite and $I=S$. It is sufficient 
to prove that $\l$ is strongly left cellular. 
First, (LC3) holds by Theorem~\ref{theo:lusztig}.

\medskip

Let $x$ and $y$ be two elements of $W$ such that $x \sim_L y$. Let $\G$ (respectively $C$) denote 
the two-sided (respectively left) cell containing $x$ and $y$. 
Then there exists $x=x_0$, $x_1$,\dots, $x_m=y=y_0$, 
$y_1$,\dots, $y_n=x$ in $W$ and elements $h_1$,\dots, $h_m$, $h_1'$,\dots, $h_n'$ of $\HC$ 
such that $C_{x_i}$ (respectively $C_{y_j}$) appears with a non-zero coefficient 
in the expression of $h_iC_{x_{i-1}}$ (respectively $h_j'C_{y_{j-1}}$) 
in the Kazhdan-Lusztig basis for $1 \le i \le m$ (respectively $1 \le j \le n$). 
Therefore, $y=x_m \lel \cdots \lel x_2 \lel x_1 = x = y_n' \lel \cdots \lel y_2' \lel y_1'=y$ and 
so $x_i$, $y_j \in C$. Hence, if we write
$$h_i C_{x_{i-1}} \equiv \sum_{u \in \G} \b_{i,u} C_u \mod \HC^{<_{LR} \G},$$
then $\b_{i,x_i} \neq 0$ and
$$v^{\alpb(\G)}h_i C_{x_{i-1}} T_{w_0} \equiv \sum_{u \in \G} v^{\alpb(\G)}\b_{i,u} C_u T_{w_0} \mod \HC^{<_{LR} \G},$$
Therefore, by Theorem~\ref{theo:lusztig}, 
$$\eta_{x_{i-1}} h_i 
C_{\l(x_{i-1})} \equiv \sum_{u \in \G} \eta_u \b_{i,u} C_{\l(u)}\mod \HC^{<_{LR} \G},$$
and so $\l(x_i) \lel \l(x_{i-1})$. This shows that $\l(y) \lel \l(x)$ and we can prove 
similarly that $\l(x) \lel \l(y)$. Therefore, $\l(C)$ is contained in a unique left cell $C'$. 
But, similarly, $\l(C')$ is contained in a unique left cell, and contains $C$. 
So $\l(C)=C'$ is a left cell. This shows (LC1).

\medskip

Finally the map $(\l,\eta)_C : \HC^L[C] \to \HC^L[\l(C)]$, $c_w^L \mapsto \eta_w c_{\l(w)}^L$ 
is obtained through the {\it right} multiplication 
by $v^{\alpb(\G)} T_{w_0}$. Since this right multiplication commutes with the left action of $\HC$, 
this implies (LC2).
\end{proof}

\bigskip

\begin{coro}\label{coro:left-cellular}
Let $I \in \PC_\frm(S)$ be such that {\bfit Lusztig's Conjectures~P1, P4,~P8 and~P9 hold} 
for the triple $(W_I,I,\ph_I)$. Then the pair $(\l_I^L,\eta_L^I)$ (respectively $(\r_I^R,\eta_R^I)$) 
is strongly left (respectively right) cellular.
\end{coro}

\bigskip

\begin{proof}
This follows from Theorems~\ref{theo:induction-left-cellular} and~\ref{theo:left-cellular}.
\end{proof}

\bigskip

It must be noticed that the maps $\l_I^L$ and $\r_I^R$ depend on the weight function $\ph$, even if 
it is not clear from the notation. The canonicity of their construction shows that, if $\s : W \to W$ 
is an automorphism such that $\s(S)=S$ and $\ph \circ \s = \ph$, then
\equat\label{eq:auto-cellular}
\s \circ \l_I^L = \l_{\s(I)}^L \circ \s\qquad
\text{and}\qquad 
\s \circ \r_I^R = \r_{\s(I)}^R \circ \s.
\endequat
For instance, {\it if $W$ is finite}, then $\o_0 : W \to W$ satisfies the above properties and so
\equat\label{eq:auto-cellular-omega}
\o_0 \circ \l_I^L = \l_{\o_0(I)}^L \circ \o_0\qquad
\text{and}\qquad 
\o_0 \circ \r_I^R = \r_{\o_0(I)}^R \circ \o_0.
\endequat

\def\gauche{{\mathrm{left}}}
\def\droite{{\mathrm{right}}}

\bigskip

\begin{coro}\label{coro:caracterisation}
Let $I \in \PC_\frm(S)$ be such that {\bfit Lusztig's Conjectures P1, P4,~P8 and P9 hold} 
for the triple $(W_I,I,\ph_I)$ and let $w \in W$. Then 
$$\eta^I_{R,w} v^{\alpb_{I,R}(w)} T_{w_I} C_w \equiv  
C_{\r_I^R(w)} \mod \HC_I^{<_R \o_{I,R}(w)}\HC$$
$$\eta^I_{L,w} v^{\alpb_{I,L}(w)} C_w T_{w_I} \equiv 
C_{\l_I^L(w)} \mod \HC\HC_I^{<_L \o_{I,L}(w)}.\leqno{\text{\it and}}$$
\end{coro}

\bigskip

\begin{proof}
It is sufficient to prove the second congruence. Let $b \in X_I$ and $y \in W_I$ be such that $w=by$ 
(so that $y=\proj_L^I(w)$). 
By Theorem~\ref{theo:geck}, 
$$C_{by} \equiv 
\sum_{\substack{(a,x) \in X_I \times W_I\\ \text{such that $a \le b$} \\ \text{and $x \sim_L y$}}}
p_{a,x,b,y}^I T_a C_x \mod \HC\HC_I^{<_L y}.$$
If $x \sim_L y$, then $\ab_I(x)=\ab_I(y)$ and $\ab_I(w_Ix)=\ab_I(w_Iy)$ by P4, so it follows from 
Theorem~\ref{theo:lusztig} that 
$$\eta^I_y v^{\alpb_I(y)}  C_{by} T_{w_I} 
\equiv 
\sum_{\substack{(a,x) \in X_I \times W_I\\ \text{such that $a \le b$} \\ \text{and $x \sim_L y$}}}
\eta^I_y\eta_x^I p_{a,x,b,y}^I T_a C_{\l_I(x)} \mod \HC\HC_I^{<_L y}.
$$
But, by Corollary~\ref{coro:egalite-geck}, $p_{a,x,b,y}^I=\eta_y^I\eta_x^I 
p_{a,\l_I(x),b,\l_I(y)}$, so 
$$\eta_{L,w}^I v^{\alpb_{I,L}(w)}  C_w T_{w_I} \equiv C_{\l_I^L(w)} 
\mod \HC \HC_I^{<_L y} T_{w_I}.$$
It then remains to notice that $T_{w_I}^{-1} C_x T_{w_I} = C_{\o_I(x)}$ for all $x \in W_I$, 
so that $\HC_I^{<_L y}T_{w_I} =T_{w_I} \HC_I^{<_L \o_I(y)} = \HC_I^{<_L \o_I(y)}$ 
and the result follows.
\end{proof}

\bigskip

An important consequence of the previous characterization is the following:

\bigskip

\begin{theo}\label{theo:commutation-strong}
Let $I \in \PC_\frm(S)$ be such that {\bfit Lusztig's Conjectures~P1, P4,~P8 and~P9 hold} 
for the triple $(W_I,I,\ph_I)$ and let $(\d,\mu)$ be a strongly left (respectively right) 
cellular pair. 
Then $\d \circ \r_I^R = \r_I^R \circ \d$ (respectively $\d \circ \l_I^L = \l_I^L \circ \d$). 
Moreover, $\eta_{R,\d(w)}^I=\mu_w \mu_{\r_I^R(w)} \eta_{R,w}^I$ 
(respectively $\eta_{L,w}^I=\eta_{L,\d(w)}^I$) for all $w \in W$. 
\end{theo}

\bigskip

\begin{proof}
Assume that $(\d,\mu)$ is strongly left cellular. 
Let $w \in W$. By~(LC3), we have $\d(w) \sim_R w$ and so~\cite[Theorem~1]{geck induction} 
$$\proj_I^R(\d(w)) \sim_R \proj_I^R(w).\leqno{(*)}$$
Now, let us write 
$$\eta_{R,w}^I v^{\alpb_{I,R}(w)} T_{w_I} C_w \equiv \sum_{u \sim_L w} \b_u C_u \mod \HC^{<_L w},$$
with $\b_u \in A$. Since $(\d,\m)$ is left cellular, we get
$$\mu_w \eta_{R,w}^I v^{\alpb_{I,R}(w)} T_{w_I} C_{\d(w)} 
\equiv \sum_{u \sim_L w} \b_u \mu_u C_{\d(u)} \mod \HC^{<_L \d(w)}.$$
But, by Corollary~\ref{coro:caracterisation}, we have
$$\eta_{R,w}^I v^{\alpb_{I,R}(w)} T_{w_I}  C_w \equiv C_{\r_I^R(w)} \mod \HC_I^{<_R \o_I(\proj_I^R(w))} \HC$$
and $\r_I^R(w) \sim_L w$ (because $\r_I^R$ is strongly right cellular by 
Corollary~\ref{coro:left-cellular}). Therefore, $\b_{\r_I^R(w)}=\mu_w\mu_{\r_I^R(w)}$. 
Again by Corollary~\ref{coro:caracterisation}, 
we get 
$$\eta_{R,w}^I v^{\alpb_{I,R}(w)} T_{w_I}  C_{\d(w)} \equiv \eta_{R,w}^I\eta_{R,\d(w)}^I 
C_{\r_I^R(\d(w))} \mod \HC_I^{<_R \o_I(\proj_I^R(w))} \HC$$
(by using also $(*)$). Combining these results, we get
$$C_{\r_I^R(\d(w))} - \eta_{R,w}^I\eta_{R,\d(w)}^I \mu_w \mu_{\r_I^R(w)} C_{\d(\r_I^R(w))} \in 
\mathop{\oplus}_{z \in \EC_1 \cup \EC_2 \cup \EC_3} A C_z,$$
where 
$$\EC_1=\{\d(u)~|~u \sim_L w\text{ and }u \neq \r_I^R(w)\},$$
$$\EC_2=\{u \in W~|~u <_L \d(w)\}$$
$$\EC_3=\{u \in W~|~\proj_I^R(u) <_R \o_I(\proj_I^R(w))\}\leqno{\text{and}}$$
(we have used the fact that $\HC_I^{<_R v} \HC=\oplus_{\proj_I^R(u) <_R v} ~A C_u$ 
for all $v \in W_I$: this result is due to Geck~\cite{geck induction}, see Theorem~\ref{theo:geck}). So, 
in order to prove that $\r_I^R(\d(w))=\d(\r_I^R(w))$ and $\eta_{R,\d(w)}^I=\mu_w \mu_{\r_I^R(w)} \eta_{R,w}^I$, 
we only need to show that 
$\d(\r_I^R(w)) \not\in \EC_1 \cup \EC_2 \cup \EC_3$. 

First, by definition, $\d(\r_I^R(w)) \not\in \EC_1$. Also, since $\r_I^R$ is strongly 
right cellular, we get that $\r_I^R(w) \sim_L w$ by (LC3) and so $\d(\r_I^R(w)) \sim_L \d(w)$ 
because $\d$ is left cellular (see (LC1)). So $\d(\r_I^R(w)) \not\in \EC_2$. 
Finally, $\d(\r_I^R(w)) \sim_R \r_I^R(w)$ because $\d$ is strongly left cellular 
(see (LC3)). So $\proj_I^R(\d(\r_I^R(w)) \sim_R \proj_I^R(\r_I^R(w))$ by~\cite{geck induction}. 
Since $\proj_I^R(\r_I^R(w))=\r_I(\proj_I^R(w)) \sim_{LR} \proj_I^R(w) \sim_{LR} \o_I(\proj_I^R(w))$ 
(see~\cite[Lemma~1.2]{lusztig action}) and so $\d(\r_I^R(w)) \not\in \EC_3$ by 
P4, P8 and P9.
\end{proof}

\section{Action of the cactus group}\label{sec:cactus}

\medskip

We recall here the definition of the {\it cactus group} $\cactus_W$ associated with $W$. The group 
$\cactus_W$ is the group with the following presentation:
\begin{itemize}
\item[$\bullet$] Generators: $(\t_I)_{I \in \PC_\irr(S)}$;

\item[$\bullet$] Relations: for all $I$, $J \in \PC_\irr(S)$, we have:
$$\begin{cases}
(\Crm 1)\quad \t_I^2=1, & \\
(\Crm 2)\quad [\t_I,\t_J]=1 & \text{if $W_{I \cup J}=W_I \times W_J$,} \\
(\Crm 3)\quad \t_I\t_J=\t_J\t_{\o_J(I)} & \text{if $I \subset J$.}
\end{cases}$$
\end{itemize}
By construction, the map $\t_I \mapsto w_I$ extends to a surjective morphism of groups 
$\cactus_W \surto W$ which will not be used in this paper.
The main result of this paper is the following:

\bigskip

\begin{theo}\label{theo:main}
Let $I$, $J \in \PC_\irr(S)$ be such that {\bfit $(\Hrm_I)$ and $(\Hrm_J)$ hold}. Then:
\begin{itemize}
\itemth{a} $[\l_I^L,\r_J^R]=\Id_W$.

\itemth{b} $(\l_I^L)^2=(\r_I^R)^2=\Id_W$.

\itemth{c} If $W_{I \cup J}=W_I \times W_J$, then $[\l_I^L,\l_J^L]=[\r_I^R,\r_J^R]=\Id_W$.

\itemth{d} If $I \subset J$, then $\l_I^L\l_J^L=\l_J^L\l_{\o_J(I)}^L$ and 
$\r_I^R\r_J^R=\r_J^R\r_{\o_J(I)}^R$.
\end{itemize}
\end{theo}

\medskip

\begin{proof}
(a) follows from Theorem~\ref{theo:commutation-strong}, while (b) is obvious. 

\smallskip

(c) Assume that $W_{I \cup J}=W_I \times W_J$. We only need to prove 
that $[\l_I^L,\l_J^L]=\Id_W$, the proof of the other equality being similar. 
Let $w \in W$ and write $w=x w'$, with $x \in X_{I \cup J}$ and 
$w' \in W_{I \cup J}$. Since $W_{I \cup J} = W_I \times W_J$ and so there exists 
$w_1 \in W_I$ and $w_2 \in W_J$ such that $w'=w_1w_2=w_2w_1$. Note also that 
$xw_1 \in X_J$, $x\l_I(w_1) \in X_J$, $xw_2 \in X_I$ and $x\l_J(w_2) \in X_I$. Therefore,
$$\l_I^L(\l_J^L(w))=\l_I^L(xw_1\l_J(w_2))=\l_I^L(x\l_J(w_2)w_1)=x\l_J(w_2)\l_I(w_1)$$
and, similarly,
$$\l_J^L(\l_I^L(w))=x\l_I(w_1)\l_J(w_2).$$
So $[\l_I^L,\l_J^L]=\Id_W$, as desired. 

\smallskip

(d) Assume here that $I \subset J$. It is easily checked that we may assume that $W$ is finite 
and $J=S$. Let $w \in W$. Then
\eqna
\l_S^L(\l_I^L(w)) &=& \r_S^L(\o_0(\l_I^L(w))) \qquad\text{by Theorem~\ref{theo:lusztig}.}\\
&=& \r_S \l_{\o_0(I)}^L(\o_0(w)) \qquad\text{by~(\ref{eq:auto-cellular-omega})},\\
&=& \l_{\o_0(I)}^L(\r_S(\o_0(w)) \qquad\text{by (a)},\\
&=& \l_{\o_0(I)}^L(\l_S(w))\qquad\text{by Theorem~\ref{theo:lusztig}.}
\endeqna
This proves the first equality and the second follows from a similar argument.
\end{proof}

\bigskip

Let $\SG_W$ denote the symmetric group on the set $W$ and assume until the end of this section that 
{\it Lusztig's Conjectures P1,~P4,~P8 and~P9 hold for the triple $(W_I,I,\ph_I)$ for all $I \in \PC_\frm(S)$}. 
The statements (b), (c) and (d) of the previous Theorem~\ref{theo:main} show that 
there exists a unique morphism of groups 
$$\fonctio{\cactus_W}{\SG_W}{\t}{\t_\ph^L}$$
such that 
$$\t_{I,\ph}^L=\l_I^L$$
for all $I \in \PC_\irr(S)$. Note that we have here emphasized the fact that the map 
depends on $\ph$. The same statements also show that 
there exists a unique morphism of groups 
$$\fonctio{\cactus_W}{\SG_W}{\t}{\t_\ph^R}$$
such that 
$$\t_{I,\ph}^R=\r_I^R$$ 
for all $I \in \PC_\irr(S)$. Moreover, Theorem~\ref{theo:main}(a) 
shows that both actions commute or, in other words, that the map
\equat\label{eq:cactus-double}
\fonctio{\cactus_W \times \cactus_W}{\SG_W}{(\t_1,\t_2)}{\t_{1,\ph}^L\t_{2,\ph}^R}
\endequat
is a morphism of groups. Let us summarize the properties of this morphism which are proved 
in this paper:

\bigskip

\begin{theo}\label{theo:action cactus}
Assume that {\bfit Lusztig's Conjectures P1,~P4,~P8 and~P9 hold for the triple $(W_I,I,\ph_I)$ 
for all $I \in \PC_\frm(S)$}. 
Let $\t \in \cactus_W$. Then there exist two sign maps $\eta_L^{\t,\ph} : W \to \mub_2$ 
and $\eta_R^{\t,\ph} : W \to \mub_2$ such that the pairs $(\t_\ph^L,\eta_L^{\t,\ph})$ 
and $(\t_\ph^R,\eta_R^{\t,\ph})$ are respectively strongly left cellular and 
strongly right cellular. 

Moreover, if $\t' \in \cactus_W$, then $[\t_\ph^L,\t_\ph^{\prime R}]=\Id_W$.
\end{theo}

\bigskip

Note that we do not claim that the sign maps in the above theorem are unique. They are obtained 
by decomposing $\t$ as a product of the generators and then compose the cellular pairs 
according to this decomposition: the resulting sign map might depend on the 
chosen decomposition. It must be added that the maps $\t_\ph^L$ and $\t_\ph^R$ depend heavily on $\ph$ 
(see for instance the case where $|S|=2$ in Section~\ref{sec:diedral}).

\bigskip

\begin{coro}\label{coro:losev}
If $W$ is a finite Weyl group and $\ph$ is constant, then the above action of $\cactus_W \times \cactus_W$ 
coincides with the one constructed by Losev~\cite[Theorem~1.1]{losev}.
\end{coro}

\bigskip

\begin{proof}
This follows from~\cite[Theorem~1.1~and~Lemma~4.7]{losev}.
\end{proof}

\section{The example of dihedral groups}\label{sec:diedral}

\medskip

\boitegrise{{\bf Hypothesis.} {\it In this section, and only in this section, 
we assume that $|S|=2$ and we write $S=\{s,t\}$. We denote by $m$ the order of 
$st$ and we assume that $3 \le m < \infty$. We denote by $\s_{s,t} : W \to W$ 
the unique involutive 
automorphism of $W$ which exchanges $s$ and $t$.}}{0.75\textwidth}

\bigskip

Recall~\cite[Proposition~5.1]{geck diedral} that Lusztig's Conjectures P1, P2,\dots, P15 
hold in this case, so that the maps $\l$ and $\r$ are well-defined. 
We aim to compute explicitly the maps $\l$ and $\r$. 
As we will see, the maps $\l$ and $\r$ depend on the weight function $\ph$. We will also 
compute the sign map $\eta$ and get the following result:

\bigskip

\begin{prop}\label{prop:sign-diedral}
If $|S|=2$ and $W$ is finite, then the sign map $\eta$ is constant on two-sided cells.
\end{prop}

\bigskip

We will need the following notation:
$$\G=W \setminus \{1,w_0\},\quad \G_s = \{w \in \G~|~ws < s\}\quad\text{and}\quad
\G_t=\{w \in \G~|~wt < t\}.$$
Note that $\G=\G_s \dotcup \G_t$, where $\dotcup$ means disjoint union.

\bigskip

\begin{rema}\label{rem:duflo-diedral}
Let
$$\DC=\{w \in W~|~\ab(w)=-\val(p_{1,w}^*)\}.$$
From P13, there exists a unique map 
$$\db : W \longto \DC$$
such that $w \sim_L \db_w$ for all $w \in W$. Its fibers are the left cells. 
Finally, it follows from~\cite[\S{2.6}]{lusztig action} 
that 
\equat\label{eq:duflo}
\r(d) =w_0 \db_{w_0 d}\qquad\text{and}\qquad \l(d)=\db_{w_0 d} w_0.
\endequat
for all $d \in \DC$.\finl
\end{rema} 

\bigskip

We define inductively two sequences 
$(s_i)_{i \ge 0}$ and $(t_i)_{i \ge 0}$ as follows:
$$
\begin{cases}
s_0=t_0=1, \\
s_{i+1}=t_i s\quad \text{and}\quad t_{i+1}=s_i t, & \text{if $i \ge 0$.}\\
\end{cases}
$$
Note that $s_1=s$, $t_1=t$ and $s_m=t_m=w_0$. Then
\equat\label{eq:gamma-s}
\G_s=\{s_1,s_2,\dots,s_{m-1}\}\qquad\text{and}\qquad \G_t=\{t_1,t_2,\dots,t_{m-1}\}.
\endequat

\bigskip

\subsection{The equal parameter case}
We assume here, and only here, that $\ph(s)=\ph(t)$ and we may also assume that $\AC=\ZM$ 
and $\ph(s)=\ph(t)=1$ (see for instance~\cite[Proposition~2.2]{bonnafe semicontinu}). 
Then~\cite[\S{8.7}]{lusztig hecke} the two-sided cells of $W$ are
$$\{1\},\quad \G\quad \text{and} \quad\{w_0\}$$
while the left cells are 
$$\{1\},\quad \G_s,\quad \G_t\quad \text{and} \quad\{w_0\}.$$
Note that $w_0\G=\G w_0=\G$. 

\bigskip

\begin{prop}\label{prop:diedral-egal}
Assume that $\ph$ is constant. Then
$$
\begin{cases}
\l(w)=\r(w)=w, & \text{if $w \in \{1,w_0\}$}\\
\l(w)=\s_{s,t}(w)w_0 & \text{if $w \not\in \{1,w_0\}$.}\\
\r(w)=w_0\s_{s,t}(w) & \text{if $w \not\in \{1,w_0\}$.}\\
\end{cases}
$$
Moreover, 
$$\eta_w=
\begin{cases}
(-1)^m & \text{if $w=1$,}\\
1 & \text{if $w =w_0$,}\\
-1 & \text{if $w \not\in \{1,w_0\}$,}\\
\end{cases}
$$
\end{prop}

\bigskip

\begin{rema}\label{rem:concret-egal}
More concretely, the (non-trivial parts of) the maps $\l$ and $\r$ are given as follows. 
If $1 \le i \le m-1$, then:
\begin{itemize}
\itemth{a} $\r(s_i)=s_{m-i}$ and $\r(t_i)=t_{m-i}$.

\itemth{b} If $m$ is even, then $\l=\r$.

\itemth{b'} If $m$ is odd, then $\l(s_i)=t_{m-i}$ and $\l(t_i)=s_{m-i}$.
\end{itemize}
In particular, if $m$ is even, then $\l$ stabilizes all the left cells 
(but nevertheless induces a non-trivial left cellular map) while, if $m$ is odd, then 
$\l$ exchanges the left cells $\G_s$ and $\G_t$ (and stabilizes all the others).\finl
\end{rema}

\bigskip

\begin{proof}
By Theorem~\ref{theo:lusztig}, we only need to compute $\r$. 
It follows from Example~\ref{ex:1} 
that $\l(1)=\r(1)=1$, that $\l(w_0)=\r(w_0)=w_0$ and that $\eta_1=(-1)^m$ and $\eta_{w_0}=1$. 
Let us also recall the following result from~\cite[\S{7}]{lusztig hecke}:
$$C_{t_i} C_s =
\begin{cases}
C_{s_2} & \text{if $i = 1$,}\\
C_{s_{i+1}} + C_{s_{i-1}} & \text{if $2 \le i \le m-1$.}\\
\end{cases}
\leqno{(*)}$$
We will use $(*)$ to show by induction on 
$i$ that $\r(s_i)=s_{m-i}$, that $\r(t_i)=t_{m-i}$ and that $\eta_{s_i}=\eta_{t_i}=-1$ 
(for $1 \le i \le m-1$). Let us first prove it for $i=1$. 
Note that $\DC \cap \G_s=\{s\}$ and $\DC \cap \G_t=\{t\}$ (see~\cite[\S{8.7}]{lusztig hecke}). 
So it follows from~(\ref{eq:duflo}) 
that $\r(s)=w_0\db_{w_0 s} $. But $w_0s \in \G_t$, so $\db_{w_0s}=t=\s_{s,t}(s)$. Therefore, 
$\r(s_1)=\r(s)=w_0\s_{s,t}(s)=s_{m-1}$, as desired. Applying the automorphism $\s_{s,t}$, 
we get $\r(t_1)=t_{m-1}$.
Note also that $\ab(w)=\ab(w_0w)$ for all $w \in \G$ (because $w_0\G=\G$). Moreover, 
$\eta_s=\eta_t=-1$ by~\cite[Theorem~2.5]{lusztig action}. 

Now, by using $(*)$, we get $C_{s_2}=C_{ts}=C_t C_s$ and 
\eqna
T_{w_0} C_{s_2}=T_{w_0} C_t C_s &\equiv& -C_{t_{m-1}} C_s \mod \HC^{<_{LR} \G} \\
&\equiv& -C_{s_{m-2}} \mod \HC^{<_{LR} \G}.
\endeqna
So $\r(s_2)=s_{m-2}$ and $\eta_{s_2}=\eta_{s_1}=-1$. Applying the automorphism $\s_{s,t}$, 
we get $\r(t_2)=t_{m-2}$ and $\eta_{t_2}=-1$.

Now, assume that $2 \le i \le m-2$ and that $\r(s_i)=s_{m-i}$, that $\r(t_i)=t_{m-i}$ 
and that $\eta_{s_i}=\eta_{t_i}=-1$. Then, by using $(*)$, we get 
\eqna
T_{w_0} C_{s_{i+1}} = T_{w_0}(C_{t_i} C_s - C_{s_{i-1}}) 
&\equiv& - C_{t_{m-i}} C_s + C_{s_{m+1-i}} \mod \HC^{<_{LR} \G} \\
&\equiv& - C_{s_{m-1-i}} \mod \HC^{<_{LR} \G}.
\endeqna
So $\r(s_{i+1})=s_{m-1-i}$ and $\eta_{s_{i+1}}=\eta_{s_i}=-1$. Applying the automorphism $\s_{s,t}$, 
we get $\r(t_{i+1})=t_{m-1-i}$ and $\eta_{t_{i+1}}=-1$. This completes the computation of $\r$. 
\end{proof}

\bigskip

\begin{rema}\label{rem:star-operation}
Note that the left cellular map $\l$ obtained here is exactly the left cellular map $w \mapsto \wti$ 
defined by Lusztig~\cite[\S{10}]{lusztig affine}. If $m=3$, this is the {\it $*$-operation}
defined by Kazhdan and Lusztig~\cite{KL}. See also~\cite[Remark~4.3~and~Example~6.3]{bonnafe geck}.\finl
\end{rema}

\bigskip

\subsection{The unequal parameter case}
Assume here, and only here, that $\ph(s) < \ph(t)$. Note that this forces $m$ to be even 
(and $m \ge 4$). 
We write $a=\ph(s)$ and $b=\ph(t)$. 
We set 
$$\G_s^<=\G_s \setminus \{s\},\quad\G_t^< = \G_t \setminus\{w_0s\}
\quad\text{and}\quad \G^<=\G_s^< \dotcup \G_t^<.$$
Then~\cite[\S{8.8}]{lusztig hecke} the two-sided cells of $W$ are
$$\{1\},\quad \{s\},\quad\G^<,\quad\{w_0s\}\quad \text{and} \quad\{w_0\}.$$
The left cells are 
$$\{1\},\quad\{s\},\quad \G_s^<,\quad \G_t^<,\quad\{w_0s\}\quad \text{and} \quad\{w_0\}.$$
Note that
$$\G_s^<=\{s_2,s_3,\dots,s_{m-1}\}\qquad\text{and}\qquad\G_t^<=\{t_1,t_2,\dots,t_{m-2}\}.$$

\bigskip

\begin{prop}\label{prop:diedral-inegal}
Assume that $\ph(s) < \ph(t)$. Let $m'=m/2$. Then
$$
\begin{cases}
\l(w)=\r(w)=w, & \text{if $w \in \{1,s,w_0s,w_0\}$,}\\
\l(s_{2i})=\r(s_{2i})= s_{m-2i} & \text{if $1 \le i \le m'-1$,}\\
\l(s_{2i+1})=\r(s_{2i+1})= s_{m+1-2i} & \text{if $1 \le i \le m'-1$,}\\
\l(t_{2i})=\r(t_{2i})= t_{m-2i} & \text{if $1 \le i \le m'-1$,}\\
\l(t_{2i-1})=\r(t_{2i-1})= t_{m-1-2i} & \text{if $1 \le i \le m'-1$.}\\
\end{cases}
$$
Moreover, 
$$\eta_w=
\begin{cases}
1 & \text{if $w \in \{1,w_0\}$,}\\
(-1)^{m'} & \text{if $w \in \{s,w_0s\}$,}\\
-1 & \text{if $w \not\in \{1,s,w_0s,w_0\}$.}\\
\end{cases}
$$
\end{prop}

\bigskip

\begin{proof}
First, note that $\l=\r$ because $w_0$ is central in $W$. 
The facts that $\l(w)=w$ if $w \in \{1,s,w_0s,w_0\}$, 
that $\eta_1=\eta_{w_0}=1$ and that $\eta_s=\eta_{w_0s}=(-1)^{m'}$ are obvious. 
Also, let $\d_{s,t}^< : W \to W$ be the map defined by 
$$\d_{s,t}^<(w)=
\begin{cases}
w & \text{if $w \in \{1,s,w_0s,w_0\}$,}\\
ws & \text{if $w \not\in \{1,s,w_0s,w_0\}$.}
\end{cases}$$
Then $\d_{s,t}^<$ is strongly left cellular~\cite[Example~6.5]{bonnafe geck} so, by 
Theorem~\ref{theo:commutation-strong}, it commutes with $\r$. In other words,
$$\forall~w \in \G^<,~\r(ws)=\r(w)s.\leqno{(*)}$$
Recall 
from~\cite[Proposition~7.6~and~\S{8.8}]{lusztig hecke} that $\DC\cap \G_s^<=\{s_3\}$ and $\DC \cap \G_t^< = \{t\}$. 
Note also that $\ab(w)=\ab(w_0w)$ for all $w \in \G^<$ (because $\G^<=w_0\G^<$). 
It follows from~(\ref{eq:duflo}) that $\r(t)=w_0 s_3=t_{m-3}$, 
and it follows from~\cite[Theorem~2.5]{lusztig action} that $\eta_t=-1$. 
Similarly, $\r(s_3)=w_0 t= s_{m-1}$ and $\eta_{s_3}=-1$. Using $(*)$, we get that 
$\r(t_2)=\r(s_3s)=s_{m-1}s=t_{m-2}$ and $\r(s_2)=\r(t_1s)=t_{m-3}s=s_{m-2}$, as desired. 
So we have proved that
$$\r(s_2)=s_{m-2},\quad\r(s_3)=s_{m-1},\quad\r(t_1)=t_{m-3}\quad\text{and}\quad
\r(t_2)=t_{m-2}$$
and that
$$\eta_{s_2}=\eta_{s_3}=\eta_{t_1}=\eta_{t_2}=-1.$$

Now, let $\z=v^{a-b}+v^{b-a}$. It follows from~\cite[Lemma~7.5~and~Proposition~7.6]{lusztig hecke} that
$$C_{t_i} C_{st} =
\begin{cases}
C_{t_{i+2}} + \z C_{t_i} & \text{if $i \in \{1,2\}$,}\\
C_{t_{i+2}} + \z C_{t_i} + C_{t_{i-2}} & \text{if $3 \le i \le m-1$.}\\
\end{cases}
$$
Using this multiplication rule and the same induction argument as in Proposition~\ref{prop:diedral-egal}, 
we get the desired result.
\end{proof}

\bigskip
\begin{rema}\label{rem:inegal}
Assume here, and only here, that $\ph(s) > \ph(t)$. 
Using the automorphism $\s_{s,t}$ which exchanges $s$ and $t$, we deduce 
from Proposition~\ref{prop:diedral-inegal} that:
$$
\begin{cases}
\l(w)=\r(w)=w, & \text{if $w \in \{1,t,w_0t,w_0\}$,}\\
\l(s_{2i})=\r(s_{2i})= s_{m-2i} & \text{if $1 \le i \le m'-1$,}\\
\l(s_{2i-1})=\r(s_{2i-1})= s_{m-1-2i} & \text{if $1 \le i \le m'-1$,}\\
\l(t_{2i})=\r(t_{2i})= t_{m-2i} & \text{if $1 \le i \le m'-1$,}\\
\l(t_{2i+1})=\r(t_{2i+1})= t_{m+1-2i} & \text{if $1 \le i \le m'-1$.}\\
\end{cases}
$$
$$\eta_w=
\begin{cases}
1 & \text{if $w \in \{1,w_0\}$,}\\
(-1)^{m'} & \text{if $w \in \{t,w_0t\}$,}\\
-1 & \text{if $w \not\in \{1,t,w_0t,w_0\}$.}\\
\end{cases}\leqno{\text{Moreover,}}
$$
This completes the proof of Proposition~\ref{prop:sign-diedral}.\finl
\end{rema}

%

\end{document}